\documentclass[10pt]{article}

\usepackage{amssymb,amsmath,amsthm}

\usepackage{hyperref}

\newtheorem{thm}{Theorem}

\newtheorem{cor}{Corollary}
\newtheorem{remark}{Remark}
\newtheorem{defi}{Definition}[section]

\newtheorem{proposition}{Proposition}
\newtheorem{exa}{Example}

\begin{document}

\title{On the characterization of non-degenerate foliations of pseudo-Riemannian manifolds with conformally flat leaves}

\author{Alfonso Garc\'{\i}a-Parrado G\'omez-Lobo
\thanks{E-mail address: {\tt alfonso@math.uminho.pt}}\\
Centro de Matem\'atica, Universidade do Minho\\ 
4710-057 Braga, Portugal}

\maketitle

\begin{abstract}
A necessary and sufficient condition for the leaves of a {\em non-degenerate} foliation of a pseudo-Riemannian manifold 
to be conformally flat is developed. The condition mimics the classical condition of the vanishing of the Weyl or Cotton tensor  
establishing the conformal flatness of a pseudo-Riemannian manifold in the sense that it is also formulated in terms of the vanishing
of certain tensors. These tensors play the role of the Weyl or the Cotton tensors and they are defined in terms of the 
the curvature of a linear torsion-free connection (the {\em bi-conformal connection}).
\end{abstract}


\section{Introduction}
A classical result in differential geometry states that a pseudo-Riemannian manifold of dimension higher than 
three is conformally flat if and only if the Weyl tensor vanishes. 
If the dimension is 3 then the Weyl tensor vanishes for any pseudo-Riemannian metric, be it conformally flat or not
and one must resort to the {\em Cotton tensor} to obtain a similar criterion as for the higher dimensional one. Both the Weyl and the 
Cotton tensors are defined exclusively in terms of the pseudo-Riemannian metric and hence they enable us to test algorithmically whether a given 
pseudo-Riemannian manifold is conformally flat or not.

One would like to have similar criteria to the ones just pointed out if we work with a foliation of a pseudo-Riemannian
manifold. The question here is whether, given such a foliation, we can answer in an algorithmically way if the induced metric on the leaves of the 
foliation is conformally flat. Of course it is clear that we can address this problem by an explicit computation of the metric induced 
on each leaf and then using the classical results explained in the previous paragraph. However, it may happen that the foliation is 
given in such a way that an explicit computation of the induced metric is very difficult or impossible. For example the foliation can be given as an 
{\em involutive distribution} or equivalently as a set of {\em integrable} differential forms. 
In this case to check whether the leaves of the foliation are conformally flat by the application of the classical result cannot be done in
a direct way. 

In this paper we provide an algorithmic criterion to check whether a foliation has conformally flat leaves which only requires the knowledge
of the distribution generating the foliation. We are able 
to construct a tensor which vanishes if and only if the leaves of the given foliation are conformally flat. Our criterion is thus similar to the 
classical criterion of the Weyl (Cotton) tensor vanishing to characterize conformally flat pseudo-Riemannian ma\-ni\-folds. 
The difference is that our tensor is constructed in terms a pair of orthogonal and complementary 
{\em projection morphisms} rather than the pseudo-Riemannian metric. These projection morphisms are defined in such a way that
the range of one of the projections coincide with the vector space spanned by the distribution at each point
of the manifold. This ensures that only the knowledge of the distribution is needed to apply our criterion. 
The only current limitation of our criterion is that it is only valid for a non-degenerate foliation 
(the induced metric on the leaves is non-degenerate).

The study of foliations of pseudo-Riemannian manifolds with conformally flat leaves has relevance in certain applications. For example 
if one knows in advance that a foliation with conformally flat leaves exists one may be able to find adapted
coordinated systems exploiting this.
Related to this is the existence problem of this type of foliations. This existence problem has been addressed in very particular cases
and with very specific motivations \cite{GARAT-PRICE,VALIENTE-COTTON}. The results presented in this paper enables 
us to address the existence problem in full
generality.

This paper is structured as follows: in section \ref{sec:distributions} we present some ge\-ne\-ralities about non-degenerate smooth 
distributions which are essential for this work. Section \ref{sec:bi-conformal} introduces the {\em bi-conformal connection} which is 
a linear torsionless connection which is defined in any pseudo-Riemannian manifold possessing a distribution. The main result 
of the paper is Theorem \ref{theo:cf-flat-leaves} of section \ref{sec:cf-foliation} where we present a necessary and sufficient condition 
guaranteeing that the leaves 
of a foliation are conformally flat. As stated above this condition takes the form of the vanishing of a certain tensor which play 
the role of either the Weyl or the Cotton tensor. This tensor is indeed defined by means of the curvature of the bi-conformal connection and 
it is of rank four or three depending on the dimension of the leaves. Specific applications of Theorem \ref{theo:cf-flat-leaves} are
presented in section \ref{sec:applications}. All the computations in this paper were carried out with the {\em Mathematica} 
suite {\em xAct} (see \cite{XACT,XPERM}).

\section{Non-degenerate smooth distributions on pseudo-Riemannian ma\-ni\-folds}
\label{sec:distributions}
Let $(V,\boldsymbol{g})$ be a $N$-dimensional smooth pseudo-Riemannian manifold and denote by $\nabla$ its torsion-free 
Levi-Civita connection. In this work we will mostly use abstract indices but sometimes index-free notation will be employed.
In this case boldface symbols will be used to denote tensorial quantities. Abstract indices of tensors and tensor fields 
arising from the tensor bundle of $V$ will be denoted with small Latin letters. 
Indices enclosed in round (square) brackets denote the operation of symmetrization (antisymmetrization).
The Levi-Civita connection $\nabla$ defines the 
Riemann and Ricci tensors in the standard way. Our conventions to define these quantities are laid by the relations
$$
\nabla_a\nabla_b X^c-\nabla_b\nabla_a X^c=R_{bad}^{\phantom{bad}c}X^d\;,\quad
R_{ac}=R_{abc}^{\phantom{abc}b}.
$$  
Indices are lowered and raised with the metric $g_{ab}$ and its inverse $g^{ab}$ (in index-free notation this is represented 
by $\boldsymbol{g}^\sharp$).
Assume now that we have defined a {\em non-degenerate smooth distribution} on the manifold $V$. 
Recall that a distribution is a smooth map $D:V\rightarrow T(V)$ such that for any point $x\in V$ its image $D_x$ under $D$
is a vector subspace of $T_x(V)$. The smoothness ensures that the dimension of $D_x$ does not depend on the point $x$ and shall be called 
the dimension of the distribution $D$.
The distribution $D$ is non-degenerate if the bi-linear form ${\boldsymbol g}|_x$ restricted to the vector 
space $D_x$ is non-degenerate. In this case we can decompose the vector space $T_x(V)$ in the form
\begin{equation}
T_x(V)=D_x\oplus D_x^{\perp}\;, 
\label{eq:tangentspace-decomposition}
\end{equation}
where $D_x^{\perp}$ is the orthogonal complement of $D_x$ with respect to the scalar product defined by ${\boldsymbol g}_x$.
The vector spaces $D_x$, $D_x^{\perp}$ can be respectively cha\-rac\-terized by means of a pair of orthogonal projectors ${\boldsymbol P}_x$, 
${\boldsymbol \Pi}_x$. These are endomorphisms of $T_x(V)$ (projection morphisms) and one has
\begin{eqnarray}
&& D_x={\boldsymbol P}_x(T_x(V))\;,\quad 
D^{\perp}_x={\boldsymbol\Pi}_x(T_x(V))\;,\quad
{\boldsymbol\Pi}_x+{\boldsymbol P}_x=1\!\!1\;,\label{eq:PPIproperties-1}\\
&&{\boldsymbol\Pi}_x\circ{\boldsymbol P}_x={\boldsymbol P}_x\circ{\boldsymbol \Pi}_x=0\;,\quad
{\boldsymbol\Pi}_x\circ{\boldsymbol\Pi}_x={\boldsymbol\Pi}_x\;,\quad {\boldsymbol P}_x\circ{\boldsymbol P}_x={\boldsymbol P}_x
\label{eq:PPIproperties-2}
\end{eqnarray}
Since the distribution $D$ is smooth we can define smooth sections ${\boldsymbol P}$, ${\boldsymbol\Pi}$ on the bundle of 
1-1 tensors $T^1_1(V)$ with ${\boldsymbol P}|_x={\boldsymbol P}_x$ and ${\boldsymbol \Pi}|_x={\boldsymbol \Pi}_x$, $x\in V$. If we represent 
these sections in index notation then is clear from (\ref{eq:PPIproperties-1})-(\ref{eq:PPIproperties-2}) that they must have the properties
\begin{equation}
P_{ab}P^{b}_{\phantom{b}c}=P_{ac}\;,\quad 
\Pi_{ab}\Pi^{b}_{\phantom{b}c}=\Pi_{ac}\;,\quad
P_{ab}\Pi^{b}_{\phantom{b}c}=0\;,\quad
P_{ab}+\Pi_{ab}=g_{ab}.
\label{eq:PPimetric-properties}
\end{equation}
Reciprocally, we can start by defining smooth sections fulfilling the relations given by the previous set of equations. 
These sections define smooth distributions $D$, $\tilde{D}$ by means of the definitions
$$
D_x={\boldsymbol P}|_x(T_x(V))\;,\quad 
\tilde{D}_x={\boldsymbol \Pi}|_x(T_x(V))\;,\quad
T_x(V)=D_x\oplus \tilde{D}_x.
$$
Note that $D_x$, $\tilde{D}_x$ may or may not be degenerate. However,
if $D_x$ is non-degenerate then we have the additional relation $\tilde{D}_x=D^{\perp}_x$.
Hence we can speak either of a (non-degenerate) distribution or the smooth sections ${\boldsymbol P}$ and ${\boldsymbol\Pi}$ indifferently.
Indeed we can just work with one projector as the other one is completely determined by the relation
${\boldsymbol\Pi}+{\boldsymbol P}=1\!\!1$. We will use the symbols $D({\boldsymbol P})$, $D({\boldsymbol\Pi})$ to denote the 
smooth distributions generated respectively by ${\boldsymbol P}$, ${\boldsymbol \Pi}$. The dimensions of these distributions are given by 
$$
\mbox{dim}(D({\boldsymbol P}))=\mbox{tr}({\boldsymbol P})=P^a_{\phantom{a}a}\;,\quad
\mbox{dim}(D({\boldsymbol \Pi}))=\mbox{tr}({\boldsymbol \Pi})=\Pi^a_{\phantom{a}a}.
$$
For any distribution $D$ we define the set $\mathfrak{X}(D)$ by
\begin{equation}
\mathfrak{X}(D)\equiv\{\vec{\boldsymbol\xi}\in\mathfrak{X}(V): \vec{\boldsymbol\xi}|_x\in D_x\;,\forall x\in V \}\;,
\end{equation}
where $\mathfrak{X}(V)$ is the set of smooth vector fields on $V$. A distribution is said to be involutive if for any 
$\vec{\boldsymbol\xi}_1$, $\vec{\boldsymbol\xi}_2$ in $\mathfrak{X}(D)$, the Lie bracket $[\vec{\boldsymbol\xi}_1,\vec{\boldsymbol\xi}_2]$ is also 
in $\mathfrak{X}(D)$.

\section{The Bi-conformal connection}
\label{sec:bi-conformal}
For any pair $P_{ab}$, $\Pi_{ab}$ fulfilling the conditions of (\ref{eq:PPimetric-properties}), define the following tensors
\begin{equation}
M_{abc}\equiv\nabla_bP_{ac}+\nabla_cP_{ab}-\nabla_aP_{bc}\;,\quad 
E_a\equiv M_{acb}P^{cb}\;,\quad 
W_a\equiv-M_{acb}\Pi^{cb}. 
\label{eq:define-M-E-W}
\end{equation}
In terms of these quantities we define the following object
\begin{equation}
 L^a_{\phantom{a}bc}\equiv\frac{1}{2\mbox{tr}({\boldsymbol P})}(E_bP^a_{\ c}+E_cP^a_{\phantom{a}b})+
\frac{1}{2(N-\mbox{tr}({\boldsymbol P}))}(W_b\Pi^a_{\phantom{a}c}+W_c
\Pi^a_{\phantom{a}b})+\frac{1}{2}(P^a_{\phantom{a}d}-\Pi^a_{\phantom{a}d})M^d_{\phantom{c}bc}.
\label{eq:define-L}
\end{equation}
This is a (1,2)-tensor with the symmetries $L^a_{\phantom{a}bc}=L^a_{\phantom{a}cb}$.
In \cite{BICONF-1} this tensor was used to introduce a new symmetric connection in the manifold $(V,{\boldsymbol g})$. We recall
now this definition.
\begin{defi}
Let $(V,{\boldsymbol g})$ be a pseudo-Riemannian manifold with its standard Levi-Civita connection 
$\nabla$ and let ${\boldsymbol P}$, ${\boldsymbol\Pi}$ be any pair fulfilling (\ref{eq:PPimetric-properties}). Use them to define the quantity 
$L^a_{\phantom{a}bc}$ by means of (\ref{eq:define-L}). The bi-conformal connection is by definition the 
only linear connection such that its difference with the Levi-Civita connection results in the tensor $L^a_{\phantom{a}bc}$.
\label{def:bi-conformal-connection} 
\end{defi}
\begin{remark}\em
Note that the bi-conformal connection is naturally defined for any smooth distribution, be it degenerate or not, 
because only relations (\ref{eq:PPimetric-properties}) are required. The results of this work,
however, require the distribution to be non-degenerate (see section \ref{sec:cf-foliation}) and our analysis 
will concern mostly this particular case.
\end{remark}

Recall that the difference between two linear connections is always a ``true tensor'', even though the connections themselves 
are not. If we denote the components of the Bi-conformal connection in a given frame by 
$\bar{\gamma}^{\boldsymbol a}_{\phantom{\boldsymbol a}\boldsymbol bc}$
then the previous definition means that
\begin{equation}
\bar{\gamma}^{\boldsymbol a}_{\phantom{\boldsymbol a}\boldsymbol{bc}}\equiv\gamma^{\boldsymbol a}_{\phantom{\boldsymbol a}{\boldsymbol bc}}+
L^{\boldsymbol a}_{\phantom{\boldsymbol a}{\boldsymbol{bc}}}\;,
\label{eq:bi-conformal-connection}
\end{equation}
where $\gamma^a_{\phantom{a}bc}$ are the components of the Levi-Civita connection in the same frame (here and in the following, 
we write frame-indices 
adopting explicit numerical values in bold).

We denote the covariant derivative associated to the bi-conformal connection by $\bar\nabla$. 
Given that $L^a_{\phantom{a}bc}$ is symmetric in its last two indices we may conclude that 
the bi-conformal connection has no torsion. The Riemann and Ricci tensors of the bi conformal connection are denoted
by 
$$
\bar{R}_{abc}^{\phantom{abc}d}\;,\quad
\bar{R}_{ab}\equiv \bar{R}_{acb}^{\phantom{acb}c}
$$
It is possible to check in explicit examples that $\bar{R}_{abc}^{\phantom{abc}c}\neq 0$, hence
the bi-conformal connection is not a Levi-Civita connection in general which means that $\bar{R}_{ab}$ need not be symmetric.

To gain some geometric intuition about the bi-conformal connection it is instructive to compute its components in a
frame adapted to the distributions $D({\boldsymbol P})$, $D({\boldsymbol \Pi})$. 
Let $\{\vec{\boldsymbol e}_1,\cdots, \vec{\boldsymbol e}_N\}$ be such a frame.
We adopt the following conventions
\begin{equation}
D({\boldsymbol P})=\mbox{Span}\{{\vec{\boldsymbol e}}_{1},\cdots,\vec{{\boldsymbol e}}_p \}\;,\quad
D({\boldsymbol\Pi})=\mbox{Span}\{{\vec{\boldsymbol e}}_{1+p},\cdots,{\vec{\boldsymbol e}}_{N}\}\;,
\label{eq:adapted-frame}
\end{equation}
where $p$ is the dimension of the distribution $D({\boldsymbol P})$ (hence $p=\mbox{tr}({\boldsymbol P})$).
We make the further assumption that
$D({\boldsymbol P})$ or $D({\boldsymbol \Pi})$ is non-degenerate which implies that both distributions
are in fact non-degenerate (see eq. (\ref{eq:tangentspace-decomposition})).
Since the frame is adapted to the distributions and they are non-degenerate, we have the relations
\begin{eqnarray}
&& P_{\boldsymbol{AB}}=\boldsymbol{g}(\vec{\boldsymbol e}_{\boldsymbol A},\vec{\boldsymbol e}_{\boldsymbol B})\;,\quad
P_{\boldsymbol{\alpha\beta}}=0\;,\quad
\Pi_{\boldsymbol{\alpha\beta}}=\boldsymbol{g}(\vec{\boldsymbol e}_{\boldsymbol \alpha},\vec{\boldsymbol e}_{\boldsymbol \beta})\;,\quad
\label{eq:components-P-Pi}\\
&&\Pi_{\boldsymbol{AB}}=0\;,\quad
P_{\boldsymbol{A}{\boldsymbol\beta}}=\boldsymbol{g}(\vec{\boldsymbol e}_{\boldsymbol A},\vec{\boldsymbol e}_{\boldsymbol \beta})=
\Pi_{\boldsymbol{A}{\boldsymbol\beta}}=0,\nonumber
\end{eqnarray}
where we have adopted the convention that capital boldface Latin letters represent frame indices with numerical 
values in the range $(1,p)$ and Greek bold symbols are frame indices in the range $(p+1,N)$. Our convention for the definition of the 
components of Levi-Civita connection (and any other linear connection) in a frame is
\begin{equation}
 \nabla_{\vec{\boldsymbol e}_{\boldsymbol a}}{\vec{\boldsymbol e}_{\boldsymbol b}}=
\gamma^{\boldsymbol c}_{\phantom{\boldsymbol c}{\boldsymbol{ab}}}{\vec{\boldsymbol e}_{\boldsymbol c}}.
\label{eq:connection-convention}
\end{equation}
The Lie brackets of the elements of this frame are characterized by the structure functions 
$C^{\boldsymbol c}_{\phantom{\boldsymbol c}{\boldsymbol{ab}}}$ defined by
\begin{equation}
[\vec{\boldsymbol e}_{\boldsymbol a},\vec{\boldsymbol e}_{\boldsymbol b}]=
C^{\boldsymbol c}_{\phantom{\boldsymbol c}{\boldsymbol{ab}}}{\vec{\boldsymbol e}_{\boldsymbol c}}=
\nabla_{\vec{\boldsymbol e}_{\boldsymbol a}}\vec{\boldsymbol e}_{\boldsymbol b}-
\nabla_{\vec{\boldsymbol e}_{\boldsymbol b}}\vec{\boldsymbol e}_{\boldsymbol a}\;,
\label{eq:structure-functions}
\end{equation}
where the last equation holds because the Levi-Civita connection has no torsion (indeed it holds
for any torsion-free connection).
To compute the components of the bi-conformal connection in this frame, we just need to compute the components of the tensor
$L^a_{\phantom{a}bc}$ and then use (\ref{eq:bi-conformal-connection}). In order to do this we are 
required to compute the components of $P_{ab}$, $\Pi_{ab}$, $\nabla_a P_{bc}$ and $\nabla_a\Pi_{bc}$ 
and then use eqs (\ref{eq:define-M-E-W})-(\ref{eq:define-L}). This is achieved using the standard definition of 
covariant derivative of a tensor 
\begin{eqnarray}
&&(\nabla_{\vec{\boldsymbol e}_a}{\boldsymbol P})(\vec{\boldsymbol e}_b,\vec{\boldsymbol e}_c)=
\vec{\boldsymbol e}_a({\boldsymbol P}(\vec{\boldsymbol e}_b,\vec{\boldsymbol e}_c))
-{\boldsymbol P}(\nabla_{\vec{\boldsymbol e}_a}\vec{\boldsymbol e}_b,\vec{\boldsymbol e}_b)
-{\boldsymbol P}(\vec{\boldsymbol e}_b,\nabla_{\vec{\boldsymbol e}_a}\vec{\boldsymbol e}_b)\;,\\
&&(\nabla_{\vec{\boldsymbol e}_a}{\boldsymbol\Pi})(\vec{\boldsymbol e}_b,\vec{\boldsymbol e}_c)=
\vec{\boldsymbol e}_a({\boldsymbol\Pi}(\vec{\boldsymbol e}_b,\vec{\boldsymbol e}_c))
-{\boldsymbol\Pi}(\nabla_{\vec{\boldsymbol e}_a}\vec{\boldsymbol e}_b,\vec{\boldsymbol e}_b)
-{\boldsymbol\Pi}(\vec{\boldsymbol e}_b,\nabla_{\vec{\boldsymbol e}_a}\vec{\boldsymbol e}_b)\;,
\end{eqnarray}
and replacing the covariant derivatives of the right hand sides by means of (\ref{eq:connection-convention}). 
The differentials $\vec{\boldsymbol e}_a({\boldsymbol P}(\vec{\boldsymbol e}_b,\vec{\boldsymbol e}_c))$, 
$\vec{\boldsymbol e}_a({\boldsymbol\Pi}(\vec{\boldsymbol e}_b,\vec{\boldsymbol e}_c))$ are 
replaced using the condition of $\nabla$ being the Levi-Civita connection
\begin{equation}
\vec{\boldsymbol e}_a({\boldsymbol g}(\vec{\boldsymbol e}_b,\vec{\boldsymbol e}_c))=
{\boldsymbol g}(\nabla_{\vec{\boldsymbol e}_a}\vec{\boldsymbol e}_b,\vec{\boldsymbol e}_c)
+{\boldsymbol g}(\vec{\boldsymbol e}_b,\nabla_{\vec{\boldsymbol e}_a}\vec{\boldsymbol e}_c)\;,
\end{equation}
which via (\ref{eq:components-P-Pi}) yields
\begin{eqnarray}
&&\vec{\boldsymbol e}_{\boldsymbol A}({\boldsymbol P}(\vec{\boldsymbol e}_{\boldsymbol B},\vec{\boldsymbol e}_{\boldsymbol C}))=
{\boldsymbol g}(\nabla_{\vec{\boldsymbol e}_{\boldsymbol A}}\vec{\boldsymbol e}_{\boldsymbol B},\vec{\boldsymbol e}_{\boldsymbol C})
+{\boldsymbol g}(\vec{\boldsymbol e}_{\boldsymbol B},\nabla_{\vec{\boldsymbol e}_{\boldsymbol A}}\vec{\boldsymbol e}_{\boldsymbol C})\;,\\
&& \vec{\boldsymbol e}_{\boldsymbol \alpha}({\boldsymbol \Pi}(\vec{\boldsymbol e}_{\boldsymbol \beta},\vec{\boldsymbol e}_{\boldsymbol \gamma}))=
{\boldsymbol g}(\nabla_{\vec{\boldsymbol e}_{\boldsymbol\alpha}}\vec{\boldsymbol e}_{\boldsymbol\beta},\vec{\boldsymbol e}_{\boldsymbol\gamma})
+{\boldsymbol g}(\vec{\boldsymbol e}_{\boldsymbol \beta},\nabla_{\vec{\boldsymbol e}_{\boldsymbol \alpha}}\vec{\boldsymbol e}_{\boldsymbol \gamma}).
\end{eqnarray}
From these formulae we deduce that the components of $L^a_{\phantom{a}bc}$ in the adapted frame will only be given in terms of the components of the 
Levi-Civita connection in the same frame and the something similar will happen with the components of the bi-conformal connection. The explicit 
relations are obtained after performing a certain amount of elementary algebra and they are
\begin{subequations}
\begin{eqnarray}
&&\bar{\gamma}^{{\boldsymbol \alpha}}_{\phantom{{\boldsymbol \alpha}}{\boldsymbol\beta}{\boldsymbol\gamma}}=
{\gamma}^{{\boldsymbol \alpha}}_{\phantom{{\boldsymbol \alpha}}{\boldsymbol\beta}{\boldsymbol\gamma}}\;,\quad 
\bar{\gamma}^{{\boldsymbol A}}_{\phantom{{\boldsymbol A}}{\boldsymbol B}{\boldsymbol C}}=
{\gamma}^{{\boldsymbol A}}_{\phantom{{\boldsymbol A}}{\boldsymbol B}{\boldsymbol C}}\;,\label{eq:bi-conf-connection-1}\\
&&\bar{\gamma}^{{\boldsymbol A}}_{\phantom{{\boldsymbol A}}{\boldsymbol B}{\boldsymbol\alpha}}=\mbox{TF}({\gamma}^{{\boldsymbol A}}_{\phantom{{\boldsymbol A}}{\boldsymbol B}{\boldsymbol\alpha}})\;,\quad
\bar{\gamma}^{{\boldsymbol A}}_{\phantom{{\boldsymbol A}}{\boldsymbol\alpha}{\boldsymbol B}}={\gamma}^{{\boldsymbol A}}_{\phantom{{\boldsymbol A}}{\boldsymbol\alpha}{\boldsymbol B}}-
{\gamma}^{{\boldsymbol A}}_{\phantom{{\boldsymbol A}}{\boldsymbol B}{\boldsymbol\alpha}}+\mbox{TF}({\gamma}^{{\boldsymbol A}}_{\phantom{{\boldsymbol A}}{\boldsymbol B}{\boldsymbol\alpha}})\;,\\
&&\bar{\gamma}^{{\boldsymbol\alpha}}_{\phantom{{\boldsymbol\alpha}}{\boldsymbol\beta}{\boldsymbol A}}=\mbox{TF}({\gamma}^{{\boldsymbol\alpha}}_{\phantom{{\boldsymbol\alpha}}{\boldsymbol\beta}{\boldsymbol A}})\;,\quad
\bar{\gamma}^{{\boldsymbol\alpha}}_{\phantom{{\boldsymbol\alpha}}{\boldsymbol A}{\boldsymbol \beta}}={\gamma}^{{\boldsymbol\alpha}}_{\phantom{{\boldsymbol\alpha}}{\boldsymbol A}{\boldsymbol\beta}}-
{\gamma}^{{\boldsymbol\alpha}}_{\phantom{{\boldsymbol\alpha}}{\boldsymbol\beta}{\boldsymbol A}}+
\mbox{TF}({\gamma}^{{\boldsymbol\alpha}}_{\phantom{{\boldsymbol\alpha}}{\boldsymbol\beta}{\boldsymbol A}})\;,\\
&&\bar{\gamma}^{{\boldsymbol A}}_{\phantom{{\boldsymbol A}}{\boldsymbol\alpha}{\boldsymbol\beta}}=
\frac{1}{2}\left({\gamma}^{{\boldsymbol A}}_{\phantom{{\boldsymbol A}}{\boldsymbol\alpha}{\boldsymbol\beta}}-
{\gamma}^{{\boldsymbol A}}_{\phantom{{\boldsymbol A}}{\boldsymbol\beta}{\boldsymbol\alpha}}\right)\;,\quad
\bar{\gamma}^{{\boldsymbol\alpha}}_{\phantom{{\boldsymbol\alpha}}{\boldsymbol A}{\boldsymbol B}}=
\frac{1}{2}\left({\gamma}^{{\boldsymbol\alpha}}_{\phantom{{\boldsymbol\alpha}}{\boldsymbol A}{\boldsymbol B}}-
{\gamma}^{{\boldsymbol\alpha}}_{\phantom{{\boldsymbol\alpha}}{\boldsymbol B}{\boldsymbol A}}\right)\;,
\label{eq:bi-conf-connection-4}
\end{eqnarray}
\end{subequations}
where (we follow Einstein's summation convention for frame indices)
\begin{eqnarray}
\mbox{TF}({\gamma}^{{\boldsymbol A}}_{\phantom{{\boldsymbol A}}{\boldsymbol B}{\boldsymbol\gamma}})\equiv
\frac{1}{2}\left({\gamma}^{{\boldsymbol A}}_{\phantom{{\boldsymbol A}}{\boldsymbol B}{\boldsymbol\gamma}}+{\gamma}^{{\boldsymbol D}}_{\phantom{{\boldsymbol D}}{\boldsymbol C}{\boldsymbol\gamma}}
P^{\boldsymbol{AC}}P_{\boldsymbol{DB}}\right)-\frac{1}{p}\gamma^{\boldsymbol C}_{\phantom{\boldsymbol C}{\boldsymbol C}{\boldsymbol\gamma}}\delta_{\boldsymbol B}^{\phantom{\boldsymbol B}{\boldsymbol A}}\;,\\
\mbox{TF}({\gamma}^{{\boldsymbol\alpha}}_{\phantom{{\boldsymbol\alpha}}{\boldsymbol\beta}{\boldsymbol C}})\equiv
\frac{1}{2}\left({\gamma}^{{\boldsymbol\alpha}}_{\phantom{{\boldsymbol\alpha}}{\boldsymbol\beta}{\boldsymbol C}}+{\gamma}^{{\boldsymbol\mu}}_{\phantom{{\boldsymbol\mu}}{\boldsymbol\gamma}{\boldsymbol C}}
P^{\boldsymbol{\alpha\gamma}}P_{\boldsymbol{\mu\beta}}\right)-\frac{1}{N-p}\gamma^{\boldsymbol\gamma}_{\phantom{\boldsymbol\gamma}{\boldsymbol\gamma}{\boldsymbol C}}
\delta_{\boldsymbol\beta}^{\phantom{\boldsymbol\beta}{\boldsymbol\alpha}}\;,
\end{eqnarray}
By looking at these formulae we see that the components of the bi-conformal connection restricted to each of the 
non-degenerate distributions $D({\boldsymbol P})$, 
$D({\boldsymbol\Pi})$ are just the components of the Levi-Civita connection (see eq. (\ref{eq:bi-conf-connection-1})). 
This is similar to what happens in 
the case of an {\em adapted linear connection} to the pair of distributions $D(\boldsymbol P)$, $D(\boldsymbol\Pi)$ 
(see \cite{BEJANCU-FARRAN} for the definition of this notion).    
However, we can also check by means of these formulae that the bi-conformal connection 
is not an adapted linear connection in general. For example, 
formula (\ref{eq:bi-conf-connection-4}) tells us that the quantities 
$\bar{\gamma}^{{\boldsymbol A}}_{\phantom{{\boldsymbol A}}{\boldsymbol\alpha}{\boldsymbol\beta}}$
and $\bar{\gamma}^{{\boldsymbol\alpha}}_{\phantom{{\boldsymbol\alpha}}{\boldsymbol A}{\boldsymbol B}}$ 
are different from zero in any adapted frame if the distributions $D(\boldsymbol P)$, 
$D(\boldsymbol\Pi)$ are not {\em involutive} (see eq. (\ref{eq:structure-functions})) whereas in the case of an adapted linear connection 
these quantities vanish in a certain adapted frame regardless
of the properties of the distributions $D(\boldsymbol P)$, $D(\boldsymbol\Pi)$.

The relevance of the bi-conformal connection was established in \cite{BICONF-1,BICONF-2} where it enabled 
to obtain a local invariant characterization of the family of {\em conformally separable} pseudo-Riemannian 
manifolds. In the next section we are going to show another geometric feature of the bi-conformal connection
dealing with non-degenerate foliations with conformally flat leaves.

\section{Characterization of foliations with conformally flat leaves}
\label{sec:cf-foliation}
A smooth foliation of the manifold $(V,\boldsymbol{g})$ is a family $\{\Sigma_{\tau}\}$, $\tau\in I\subset\mathbb{R}^p$ 
of smooth embedded sub-manifolds of dimension $N-p$ such that for any $\tau\neq \tau'$, 
$\Sigma_{\tau}\cap \Sigma_{\tau'}=\varnothing$ and $V=\cup_{\tau\in I}\Sigma_{\tau}$. Each of the submanifolds $\Sigma_{\tau}$
is called a leaf of the foliation and the natural number $p$ is their co-dimension. 
The foliation is {\em non-degenerate} if the induced metric (first fundamental form) 
on each leaf is non-degenerate everywhere. Note that any non-degenerate 
foliation always generates a non-degenerate distribution
defined for any $x\in V$ by 
$$
x\mapsto T_x(\Sigma_\tau)\;, x\in\Sigma_\tau 
$$
where $T_x(\Sigma_\tau)$ is the subspace of the tangent space $T_x(V)$ corresponding to the tangent space of the leaves. 
By Frobenius theorem, the distribution so defined is involutive.
Also, associated to any foliation there is a smooth map
$t:V\rightarrow \mathbb{R}^p$ defined for any $x\in V$ by $t(x)=\tau$, where $x\in\Sigma_{\tau}$. This function will be referred to
as the level set map of the foliation.

As explained in section \ref{sec:distributions}, any non-degenerate smooth distribution is characterized by a
smooth projector ${\boldsymbol P}$ and its orthogonal complement ${\boldsymbol\Pi}$, and therefore we can 
apply these ideas to the present si\-tuation.  We choose the convention that the non-degenerate distribution generated by 
$\{\Sigma_{\tau}\}$ is $D({\boldsymbol\Pi})$ which means that dim$(D({\boldsymbol\Pi}))=N-p$ 
and tr$({\boldsymbol\Pi})=N-p$, tr$({\boldsymbol P})=p$

In this work we are interested in finding necessary and sufficient conditions 
which guarantee that the leaves of a non-degenerate foliation are conformally flat.
The following theorem supplies these conditions. 
\begin{thm}
The leaves with dimension $N-p>3$, $N>3$ of a non-degenerate foliation of a pseudo-Riemannian manifold $(V,{\boldsymbol g})$ 
are conformally flat if and only if the pair of complementary projectors $({\boldsymbol\Pi},{\boldsymbol P})$ associated to the 
distribution defined 
by the foliation fulfill the following condition
\begin{equation}
{\boldsymbol T}^{||}({\boldsymbol\Pi})=0\;, 
\label{eq:bi-conformally-flat-1}
\end{equation}
where
\begin{equation}
(T^{||}({\boldsymbol \Pi}))_{cab}^{\phantom{cab}d}\equiv \Pi^d_{\phantom{d}r}\Pi^s_{\phantom{s}c}\Pi^t_{\phantom{t}a}\Pi^q_{\phantom{q}b}T_{stq}^{\phantom{stq}r},
\label{eq:define-T-parallel}
\end{equation}
\begin{equation}
T_{bac}^{\phantom{bac}d}\equiv 2\bar R_{bac}^{\phantom{bac}d}-\frac{2}{2-N-p}
(\Pi^d_{\phantom{d}c}L^\Pi_{[ab]}+\Pi^d_{\phantom{d}[b}L^\Pi_{a]c}
+\Pi_{c[a}L^\Pi_{b]q}\Pi^{qd})\;,
\label{eq:define-T}
\end{equation}
and 
\begin{eqnarray}
&& L^\Pi_{\phantom{0}bc}\equiv
2\left(\Pi^d_{\phantom{d}r}\bar R_{bdc}^{\phantom{bdc}r}-\frac{1}{N-p}(\Pi^d_{\phantom{d}c}\Pi^r_{\phantom{r}q}
\bar R_{bdr}^{\phantom{bdr}q}+\Pi^d_{\phantom{d}b}\Pi^r_{\phantom{r}q}R_{cdr}^{\phantom{cdr}q}-
\Pi^r_{\phantom{r}q}\bar R_{cbr}^{\phantom{cbr}q})\right)+\nonumber\\
&&+\frac{\bar R^\Pi}{1-N+p}\Pi_{bc}\;,\label{eq:define-L1}\\
&&\bar R^\Pi\equiv \Pi^d_{\phantom{d}r}\bar R_{bdc}^{\phantom{bdc}r}\Pi^{cb}\;,
\end{eqnarray}
the bi-conformal connection in these expressions being that defined from $P_{ab}$ and 
$\Pi_{ab}=g_{ab}-P_{ab}$ in the way 
explained in Definition \ref{def:bi-conformal-connection}. If $N-p=3$
then eq. (\ref{eq:bi-conformally-flat-1}) is to be replaced by
\begin{equation}
{\boldsymbol B}^{||}({\boldsymbol\Pi})=0\;, 
\label{eq:bi-conformally-flat-2}
\end{equation}
where
\begin{equation}
(B^{||}({\boldsymbol\Pi}))_{abc}\equiv \Pi_a^{\phantom{a}r}\Pi_b^{\phantom{b}s}\Pi_c^{\phantom{c}q}\bar\nabla_{[r}L^\Pi_{s]q}.
\label{eq:define-B}
\end{equation}
\label{theo:cf-flat-leaves}
\end{thm}

\proof We choose a frame adapted to the non-degenerate distributions $D({\boldsymbol P})$, $D({\boldsymbol\Pi})$, in 
the way explained
in section \ref{sec:bi-conformal}, (see eq. (\ref{eq:adapted-frame})) and show that in 
that frame ${\boldsymbol T}^{||}({\boldsymbol\Pi})=0$
($N-p>3$) or ${\boldsymbol B}^{||}({\boldsymbol\Pi})$ ($N-p=3$) vanish if and only if 
the leaves of the foliation are conformally flat. 
In fact, given that $D(\boldsymbol\Pi)$ is involutive
we can find a frame such that $\vec{\boldsymbol e}_{\boldsymbol\alpha}=\partial/\partial x^{\boldsymbol\alpha}$
in certain local coordinates $(x^1,\cdots,x^{N-p},x^{N-p+1},\cdots,x^N)$, $N-p+1\leq{\boldsymbol\alpha}\leq N$.
In this adapted frame the components of the induced metric on the leaves (first fundamental form) 
are just $\Pi_{\boldsymbol{\alpha\beta}}$ (see eq. (\ref{eq:components-P-Pi})) and the quantities 
$\gamma^{\boldsymbol\alpha}_{\phantom{\boldsymbol\alpha}{\boldsymbol\mu}{\boldsymbol\nu}}$
correspond to the components of the Levi-Civita connection arising from the first fundamental form.
In addition, since the connection $\nabla$ has no torsion,  
we have the relation
\begin{eqnarray}
&& 0=\left[\frac{\partial}{\partial x^{\boldsymbol\alpha}},\frac{\partial}{\partial x^{\boldsymbol\beta}}\right]=[\vec{\boldsymbol e}_{\boldsymbol\alpha},\vec{\boldsymbol e}_{\boldsymbol\beta}]
=\nabla_{\vec{\boldsymbol e}_{\boldsymbol\alpha}}\vec{\boldsymbol e}_{\boldsymbol\beta}-
\nabla_{\vec{\boldsymbol e}_{\boldsymbol\beta}}\vec{\boldsymbol e}_{\boldsymbol\alpha}=\\
&&=(\gamma^{\boldsymbol A}_{\phantom{\boldsymbol A}\boldsymbol{\alpha\beta}}-\gamma^{\boldsymbol A}_{\phantom{\boldsymbol A}\boldsymbol{\beta\alpha}})
\vec{\boldsymbol e}_{\boldsymbol A}+(\gamma^{\boldsymbol\gamma}_{\phantom{\boldsymbol\gamma}\boldsymbol{\alpha\beta}}
-\gamma^{\boldsymbol\gamma}_{\phantom{\boldsymbol\gamma}\boldsymbol{\beta\alpha}})
\vec{\boldsymbol e}_{\boldsymbol\gamma}\;,
\end{eqnarray}
which entails
$\gamma^{\boldsymbol A}_{\phantom{\boldsymbol A}\boldsymbol{\alpha\beta}}=\gamma^{\boldsymbol A}_{\phantom{\boldsymbol A}\boldsymbol{\beta\alpha}}$.
Using this in the first condition of (\ref{eq:bi-conf-connection-4}) we conclude that 
$\bar\gamma^{\boldsymbol A}_{\phantom{\boldsymbol A}\boldsymbol{\alpha\beta}}=0$. Next we need to compute the components of the bi-conformal 
curvature ${\bar R}_{abc}^{\phantom{abc}d}$ in the chosen frame. To that end we use the standard definition of the curvature
of a linear connection in differential geometry particularized to our frame 
\begin{equation}
\bar{\boldsymbol R}(\vec{\boldsymbol e}_{\boldsymbol a},\vec{\boldsymbol e}_{\boldsymbol b})\vec{\boldsymbol e}_{\boldsymbol c}=
\bar\nabla_{\vec{\boldsymbol e}_{\boldsymbol b}}\bar\nabla_{\vec{\boldsymbol e}_{\boldsymbol a}}\vec{\boldsymbol e}_{\boldsymbol c}
-\bar\nabla_{\vec{\boldsymbol e}_{\boldsymbol a}}\bar\nabla_{\vec{\boldsymbol e}_{\boldsymbol b}}\vec{\boldsymbol e}_{\boldsymbol c}-
\bar\nabla_{[\vec{\boldsymbol e}_{\boldsymbol a},\vec{\boldsymbol e}_{\boldsymbol b}]}\vec{\boldsymbol e}_{\boldsymbol c}.
\label{eq:riemann-components}
\end{equation}
In fact, if we look at the definition of ${\boldsymbol T}^{||}({\boldsymbol\Pi})$, we realize that to compute this quantity
in the adapted frame we only need 
to compute the scalars $\bar{\boldsymbol R}_{\boldsymbol{\mu\beta\gamma}}^{\phantom{\boldsymbol{\mu\beta\gamma}}\boldsymbol\lambda}$
and $L^\Pi_{\boldsymbol{\alpha\beta}}$
which are, respectively, the only non-vanishing components of the tensors 
$(\bar R^{||}({\boldsymbol\Pi}))_{cab}^{\phantom{cab}d}\equiv\Pi^d_{\phantom{d}r}\Pi^s_{\phantom{s}c}\Pi^t_{\phantom{t}a}\Pi^q_{\phantom{q}b}\bar R_{stq}^{\phantom{stq}r}$
and $(L^{\Pi||})_{ab}\equiv\Pi^t_{\phantom{t}a}\Pi^q_{\phantom{q}b}L^\Pi_{\phantom{1}tq}$ in the adapted frame.
This is achieved by combining (\ref{eq:riemann-components}), (\ref{eq:connection-convention}) (with $\nabla$ replaced by $\bar\nabla$)  
and (\ref{eq:structure-functions}) with the result
\begin{eqnarray}
&&\bar{\boldsymbol R}_{\boldsymbol{\mu\beta\gamma}}^{\phantom{\boldsymbol{\mu\beta\gamma}}\boldsymbol\lambda}=
2\vec{\boldsymbol e}_{[\boldsymbol\beta}(\bar\gamma^{\boldsymbol\lambda}_{\phantom{\boldsymbol\lambda}{\boldsymbol\mu}]{\boldsymbol\gamma}})+
2\bar\gamma^{\boldsymbol\alpha}_{\phantom{\boldsymbol\alpha}{\boldsymbol[\mu}|{\boldsymbol\gamma}}
\bar\gamma^{\boldsymbol\lambda}_{\phantom{\boldsymbol\lambda}{|\boldsymbol\mu}]{\boldsymbol\alpha}}-
\bar\gamma^{\boldsymbol\lambda}_{\phantom{\boldsymbol\lambda}{\boldsymbol\alpha}{\boldsymbol\gamma}}
C^{\boldsymbol\alpha}_{\phantom{\boldsymbol\alpha}{\boldsymbol\mu}{\boldsymbol\beta}}+
2\bar\gamma^{\boldsymbol A}_{\phantom{\boldsymbol A}{[\boldsymbol\mu|}{\boldsymbol\gamma}}
\bar\gamma^{\boldsymbol\lambda}_{\phantom{\boldsymbol\lambda}{|\boldsymbol\beta]}{\boldsymbol A}}-\nonumber\\
&&-\bar\gamma^{\boldsymbol\lambda}_{\phantom{\boldsymbol\lambda}{\boldsymbol A}{\boldsymbol\gamma}}
C^{\boldsymbol A}_{\phantom{\boldsymbol A}{\boldsymbol\mu}{\boldsymbol\beta}}.
\label{eq:bi-conformal-curvature-comps}
\end{eqnarray}
The fact that the distribution $D(\boldsymbol\Pi)$ is involutive and our choice of 
an adapted frame implies respectively that $C^{\boldsymbol A}_{\phantom{\boldsymbol A}\boldsymbol{\mu\beta}}=0$
and $C^{\boldsymbol\alpha}_{\phantom{\boldsymbol\alpha}\boldsymbol{\mu\beta}}=0$.
Also if we use (\ref{eq:bi-conf-connection-1}) and the property 
$\bar\gamma^{\boldsymbol A}_{\phantom{\boldsymbol A}\boldsymbol{\alpha\beta}}=0$ discussed above we 
can reduce (\ref{eq:bi-conformal-curvature-comps}) to
\begin{equation}
\bar{\boldsymbol R}_{\boldsymbol{\mu\beta\gamma}}^{\phantom{\boldsymbol{\mu\beta\gamma}}\boldsymbol\lambda}=
2\vec{\boldsymbol e}_{[\boldsymbol\beta}(\gamma^{\boldsymbol\lambda}_{\phantom{\boldsymbol\lambda}{\boldsymbol\mu}]{\boldsymbol\gamma}})+
2\gamma^{\boldsymbol\alpha}_{\phantom{\boldsymbol\alpha}{\boldsymbol[\mu}|{\boldsymbol\gamma}}
\gamma^{\boldsymbol\lambda}_{\phantom{\boldsymbol\lambda}{|\boldsymbol\mu}]{\boldsymbol\alpha}}=
{\boldsymbol R}_{\boldsymbol{\mu\beta\gamma}}^{\phantom{\boldsymbol{\mu\beta\gamma}}\boldsymbol\lambda}.
\label{eq:barRiemann}
\end{equation}
From the previous formula we deduce that the scalars 
${\boldsymbol R}_{\boldsymbol{\mu\beta\gamma}}^{\phantom{\boldsymbol{\mu\beta\gamma}}\boldsymbol\lambda}$ 
can be understood as the components 
of the Riemann tensor of the induced metric (first fundamental form) on the leaves of the foliation. 
Regarding the non-vanishing components of  $L^{\Pi||}$ a direct computation using 
(\ref{eq:define-L1}) shows
\begin{equation}
L^\Pi_{\phantom{1}\boldsymbol{\beta\gamma}}=2\bar R_{\boldsymbol{\beta\alpha\gamma}}^{\phantom{\boldsymbol{\beta\alpha\gamma}}\boldsymbol{\alpha}}-
\frac{2}{N-p}\bar R_{\boldsymbol{\beta\gamma\alpha}}^{\phantom{\boldsymbol{\beta\gamma\alpha}}\boldsymbol{\alpha}}+
\frac{\bar R^\Pi}{1-N+p}\Pi_{\boldsymbol{\beta\gamma}}\;,\quad
\bar R^\Pi=\bar R_{\boldsymbol{\beta\alpha\gamma}}^{\phantom{\boldsymbol{\beta\alpha\gamma}}\boldsymbol{\alpha}}\Pi^{\boldsymbol{\beta\gamma}}\;, 
\end{equation}
which on using (\ref{eq:barRiemann}) yields
\begin{equation}
L^\Pi_{\phantom{1}\boldsymbol{\beta\gamma}}=2 R_{\boldsymbol{\beta\gamma}}+\frac{\Pi_{\boldsymbol{\beta\gamma}}R}{1-N+p}\;,\quad
R_{\boldsymbol{\beta\gamma}}\equiv R_{\boldsymbol{\beta\delta\gamma}}^{\phantom{\boldsymbol{\beta\delta\gamma}}\boldsymbol{\delta}}\;,\quad
R\equiv R_{\boldsymbol{\beta\gamma}}\Pi^{\boldsymbol{\beta\gamma}}\;,
\label{eq:L1-split}
\end{equation}
where, for the same reasons pointed out above with the Riemann tensor,
$R_{\boldsymbol{\beta\gamma}}$ coincides with the Ricci tensor computed from the first fundamental form and $R$ is its 
scalar curvature. Combining (\ref{eq:barRiemann}) and (\ref{eq:L1-split}) with (\ref{eq:define-T-parallel})-(\ref{eq:define-T}) we obtain
after some algebra
\begin{equation}
(T^{||}({\boldsymbol \Pi}))_{\boldsymbol{\beta\alpha\gamma}}^{\phantom{\boldsymbol{\beta\alpha\gamma}}{\boldsymbol\delta}}= 
2 R_{\boldsymbol{\beta\alpha\gamma}}^{\phantom{\beta\alpha\gamma}{\boldsymbol\delta}}-\frac{2}{2-N-p}
(\Pi^{\boldsymbol\delta}_{\phantom{\boldsymbol\delta}[{\boldsymbol\beta}}L^\Pi_{\boldsymbol{\alpha}]{\boldsymbol\gamma}}
+\Pi_{\boldsymbol{\gamma}[\boldsymbol{\alpha}}L^\Pi_{\boldsymbol{\beta}]\boldsymbol{\mu}}\Pi^{\boldsymbol{\mu\delta}})=
2 W_{\boldsymbol{\gamma\alpha\beta}}^{\phantom{\boldsymbol{\gamma\alpha\beta}}{\boldsymbol\delta}},
\end{equation}
where $W_{\boldsymbol{\gamma\alpha\beta}}^{\phantom{\boldsymbol{\gamma\alpha\beta}}{\boldsymbol\delta}}$
is the Weyl tensor defined by the leaves first fundamental form. Hence, when $N-p>3$ the tensor $T^{||}({\boldsymbol \Pi})$ vanish if and 
only if these leaves are conformally flat. In a similar fashion as above, from (\ref{eq:define-B}) we deduce the only non-vanishing components of 
$B^{||}({\boldsymbol P})$ are $B^{||}({\boldsymbol P})_{\boldsymbol{\alpha\beta\gamma}}$ and these are given by
\begin{equation}
B^{||}({\boldsymbol P})_{\boldsymbol{\alpha\beta\gamma}}=
\bar\nabla_{[\boldsymbol{\gamma}}L^\Pi_{\boldsymbol{\beta}]\boldsymbol{\mu}}=
\nabla_{[\boldsymbol{\gamma}}L^\Pi_{\boldsymbol{\beta}]\boldsymbol{\mu}}\;,
\end{equation}
where in the last step we used again (\ref{eq:bi-conf-connection-1}) and 
the property $\bar\gamma^{\boldsymbol A}_{\phantom{\boldsymbol A}\boldsymbol{\alpha\beta}}=0$ in the expansion of 
$\bar\nabla_{[\boldsymbol{\gamma}}L^\Pi_{\boldsymbol{\beta}]\boldsymbol{\mu}}$
in terms of the components of the bi-conformal connection. 
Since $L^\Pi_{\boldsymbol{\beta}\boldsymbol{\mu}}$
is, up to a factor, the Schouten tensor of the induced metric on the leaves 
we conclude that $\nabla_{[\boldsymbol{\gamma}}L^\Pi_{\boldsymbol{\beta}]\boldsymbol{\mu}}$
is the Cotton tensor of the same metric and therefore $B^{||}({\boldsymbol P})$ will be zero if and only if the Cotton tensor 
of these leaves vanish. In the case of the leaves being of dimension $3$ ($N-p=3$) this is equivalent to their being conformally flat.\qed

\begin{cor}Under the hypotheses of theorem \ref{theo:cf-flat-leaves}
the tensor $\bar R^{||}({\boldsymbol P})$ is zero if and only if 
the leaves of the foliation generated by the pair $({\boldsymbol P},{\boldsymbol \Pi})$ are flat.
\end{cor}
\proof See considerations coming after eq. (\ref{eq:barRiemann}). \qed

\begin{remark}\em
If $N\leq 3$ then the leaves of any foliation are of dimension 1 or 2 and hence they are trivially 
conformally flat. Therefore to deal with non-trivial situations we need to consider the case $N>3$ which is the situation assumed 
by the hypotheses of Theorem \ref{theo:cf-flat-leaves}. 
\end{remark}

\section{Applications}
\label{sec:applications}
The obvious application of our result comes in the determination of the conformal flatness of the leaves of a foliation. 
Of course it is fairly trivial to determine this if the 
foliation is given in an {\em explicit form} (for example in terms of the level set functions) because in this 
case it is possible to compute {\em explicitly} the first fundamental form for each leaf and then one just needs 
to compute the Weyl or the Cotton tensor for it. However, there are situations where one does not have a explicit form 
of the foliation or the first fundamental form of the leaves is very difficult to compute. 
For example, this situation arises if one is interested in investigating the existence 
of a non-degenerate foliation of a pseudo-Riemannian manifold by conformally flat leaves. To understand better this problem
let us consider the special simpler case where the leaves of the foliation have co-dimension $p=1$. In this case the distribution 
$D(\boldsymbol\Pi)$ can be characterized by a 1-form ${\boldsymbol\omega}\in\Lambda^1(V)$ such 
that ${\boldsymbol\omega}_x(\vec{\boldsymbol X}_x)=0$, for any $\vec{\boldsymbol X}_x\in D_x({\boldsymbol\Pi})$, and $ x\in V$.
In addition, if the distribution $D({\boldsymbol\Pi})$ is non-degenerate and involutive then, 
${\boldsymbol\omega}\wedge d{\boldsymbol\omega}=0$ (integrable 1-form) and 
${\boldsymbol g}^\sharp({\boldsymbol\omega},{\boldsymbol\omega})\neq 0$.
Indeed it is possible to establish a 
relation between ${\boldsymbol\omega}$ and the projectors ${\boldsymbol P}$, ${\boldsymbol\Pi}$.
\begin{proposition}
Any non-degenerate foliation of a pseudo-Riemannian manifold $(V,{\boldsymbol g})$
with co-dimension 1 leaves can be characterized by orthogonal projectors ${\boldsymbol P}$, 
${\boldsymbol\Pi}$ adopting the form
\begin{equation}
{\boldsymbol P}= \frac{{\boldsymbol\omega}^\sharp\otimes{\boldsymbol\omega}}
{{\boldsymbol g}^\sharp({\boldsymbol\omega},{\boldsymbol\omega})}\;,
\quad {\boldsymbol\Pi}=1\!\!1-{\boldsymbol P}\;, 
\label{eq:definePiP}
\end{equation}
where ${\boldsymbol\omega}\in\Lambda^1(V)$ is such that 
${\boldsymbol\omega}\wedge d{\boldsymbol\omega}=0$ and ${\boldsymbol\omega}^\sharp$ is the vector determined by the relation
${\boldsymbol\omega}^\sharp(\cdot)={\boldsymbol g}^{\sharp}({\boldsymbol\omega},\cdot)$, 
${\boldsymbol g}^\sharp({\boldsymbol\omega},{\boldsymbol\omega})\neq 0$.
\end{proposition}

\proof 
For any $\boldsymbol{\omega}\in\Lambda^1(V)$ with ${\boldsymbol g}^\sharp({\boldsymbol\omega},{\boldsymbol\omega})\neq 0$
the quantities ${\boldsymbol \Pi}$ and 
${\boldsymbol P}$ introduced by (\ref{eq:definePiP}), define the non-degenerate distributions
\begin{equation}
 D({\boldsymbol P})=\mbox{Span}\{{\boldsymbol\omega}^\sharp\}\;,
\end{equation}
and its orthogonal complement $D({\boldsymbol\Pi})$.
Also, by construction ${\boldsymbol\omega}(\vec{\boldsymbol X})=0$, for any $\vec{\boldsymbol X}_x\in D({\boldsymbol\Pi})_x$.
If in addition ${\boldsymbol\omega}\wedge d{\boldsymbol\omega}=0$, then $D({\boldsymbol\Pi})$ is an involutive
distribution and by Frobenius theorem it generates a foliation of $V$ with co-dimension 1 leaves.
Since $D({\boldsymbol\Pi})$ is non-degenerate, the foliation is also non-degenerate.
Reciprocally, any non-degenerate foliation of $V$ with co-dimension 1 leaves determines an integrable 1-form ${\boldsymbol\omega}$,
with ${\boldsymbol g}^\sharp({\boldsymbol\omega},{\boldsymbol\omega})\neq 0$ as explained above. 
We can then take then equation (\ref{eq:definePiP}) as the definition for ${\boldsymbol P}$ and
${\boldsymbol\Pi}$. 
\qed

Now if we wish to adapt the result of Theorem \ref{theo:cf-flat-leaves} to the case of a foliation with leaves of
co-dimension 1 then it is clear that one can write condition (\ref{eq:bi-conformally-flat-1}) or 
(\ref{eq:bi-conformally-flat-2}) as a condition about the existence of a certain 1-form $\boldsymbol\omega$. 
This condition will be generically not satisfied by 
any $\boldsymbol\omega$ unless the pseudo-Riemannian manifold admits a foliation by co-dimension 1 hypersurfaces.
See example \ref{ex:schw} for the study of a particular case of this situation.  
Another application arises if we are given an integrable non-degenerate 1-form and we are asked whether the foliation 
it generates is formed by conformally 
flat leaves. In this case we do not have an explicit form for the leaves, nor their first fundamental form
and hence the {\em classical procedure} cannot be applied directly. This problem does not arise if we 
apply Theorem \ref{theo:cf-flat-leaves} as the explicit form of the leaves is not required. With just the knowledge 
of the integrable 1-form we can build the projectors $\boldsymbol{P}$, $\boldsymbol{\Pi}$ via eqs. 
(\ref{eq:definePiP}) and then check whether 
(\ref{eq:bi-conformally-flat-1}) or (\ref{eq:bi-conformally-flat-2}) are fulfilled. 

\begin{exa}
\label{ex:schw}
\em
 Take as pseudo-Riemannian manifold $(V,{\boldsymbol g})$ the region of the Schwarzschild solution of the vacuum Einstein equations
covered by the standard Schwarzschild coordinates
$(t,r,\theta,\phi)$
\begin{equation}
 {\boldsymbol g}=-\left(1-\frac{2M}{r}\right)dt\otimes dt+\frac{1}{1-\frac{2M}{r}}dr\otimes dr+r^2
(d\theta\otimes d\theta+\sin^2\theta d\phi\otimes d\phi)\;,
\end{equation}
where $M$ is the mass parameter. Suppose we ask ourselves about the existence of a foliation of $V$ whose leaves are
conformally flat and given by the family of hypersurfaces
$$
\{\Phi(t,r)=c\}\;,\quad c\in\mathbb{R}\;,\Phi\in C^2(V).
$$
This means that the level set function for such a foliation must take 
the form $\Phi(t,\phi)$ and hence we can choose the 1-form
$$
{\boldsymbol\omega}=d\Phi=\Phi_t dt+\Phi_\phi d\phi\;,\quad \Phi_t\equiv\frac{\partial\Phi}{\partial t}\;,\quad
\Phi_\phi\equiv\frac{\partial\Phi}{\partial \phi}
$$
to construct the pair $({\boldsymbol P},{\boldsymbol\Pi})$ through eq. (\ref{eq:definePiP}). Since 
$N=\mbox{dim}(V)=4$ we need to check the condition ${\boldsymbol B}^{||}({\boldsymbol P})=0$ in accordance with Theorem 
\ref{theo:cf-flat-leaves}. Once the pair $({\boldsymbol P},{\boldsymbol\Pi})$ is given this is an algorithmic computation which 
can be carried out with the system \cite{XACT}. The result is an overdetermined system of partial differential equations for the 
function $\Phi$ which we need to study. One of these equations is 
$$
(B^{||})_{\phi\theta\phi}=\frac{4 r^5 (2M-r)\Phi_{\phi}^2 \Phi_{t}^4\left(M (9M-4 r) \Phi_{\phi}^2+r^4\Phi _{t}^2\right)
\sin^5\theta\cos\theta}
{\left((2M-r) \Phi_{\phi}^2+r^3 \Phi_{t}^2\sin^2\theta\right)^4}=0\;,
$$
which is only true if the numerator of the rational expression vanishes. As $\Phi$ is a function of just $t$ and $\phi$
the only possibilities are that either $\Phi_{\phi}=0$ or $\Phi_{t}=0$. The latter possibility can be discarded as one can 
explicitly check that there are components of ${\boldsymbol B}^{||}({\boldsymbol P})$ which cannot vanish, so we are left 
with the case
$\Phi_{\phi}=0$ which makes ${\boldsymbol B}^{||}({\boldsymbol P})$ vanish identically. This corresponds to the trivial 
case of the leaves
being the standard static hypersurfaces which, as is well-known, are conformally flat for the Schwarzschild solution. 
\end{exa}

\section*{Acknowledgments}
Financial support by the Research Centre of Mathematics of the University of
Minho (Portugal) through the ``Funda\c{c}\~ao para a Ci\^encia e a Tecnolog\'{\i}a (FCT) Pluriannual
Funding Program'' and by project CERN/FP/123609/2011 is gratefully acknowledged.

\bibliographystyle{amsplain}
\bibliography{/home/alfonso/trabajos/BibDataBase/Bibliography}

\providecommand{\bysame}{\leavevmode\hbox to3em{\hrulefill}\thinspace}
\providecommand{\MR}{\relax\ifhmode\unskip\space\fi MR }
\providecommand{\MRhref}[2]{%
  \href{http://www.ams.org/mathscinet-getitem?mr=#1}{#2}
}
\providecommand{\href}[2]{#2}
\begin{thebibliography}{1}

\bibitem{BEJANCU-FARRAN}
A.~Bejancu and H.~R. Farran, \emph{Foliations and {G}eometric {S}tructures},
  first ed., Mathematics and {I}ts {A}pplications, Springer, Dordrecht, 2006.

\bibitem{GARAT-PRICE}
A.~Garat and R.~H. Price, \emph{Nonexistence of conformally flat slices of the
  {K}err spacetime}, Phys. Rev. D (3) \textbf{61} (2000), no.~12, 124011, 4.

\bibitem{BICONF-1}
A.~Garc{\'{\i}}a-Parrado, \emph{Bi-conformal vector fields and the local
  geometric characterization of conformally separable pseudo-{R}iemannian
  manifolds. {I}}, J. Geom. Phys. \textbf{56} (2006), no.~7, 1069--1095.

\bibitem{BICONF-2}
\bysame, \emph{Bi-conformal vector fields and the local geometric
  characterization of conformally separable pseudo-{R}iemannian manifolds.
  {II}}, J. Geom. Phys. \textbf{56} (2006), no.~9, 1600--1622.

\bibitem{XACT}
J.~M. Mart\'{\i}n-Garc\'{\i}a, \emph{x{A}ct: efficient tensor computer
  algebra}, \url{http://www.xact.es}.

\bibitem{XPERM}
\bysame, \emph{x{P}erm: fast index canonicalization for tensor computer
  algebra}, Computer Physics Communications \textbf{179} (2008), 597--603.

\bibitem{VALIENTE-COTTON}
J.~A. Valiente~Kroon, \emph{Asymptotic expansions of the {C}otton-{Y}ork tensor
  on slices of stationary spacetimes}, Classical Quantum Gravity \textbf{21}
  (2004), no.~13, 3237--3249.

\end{thebibliography}

\end{document}